\documentclass[reqno]{amsart}

\usepackage{amsmath,amssymb,amsthm}

\usepackage{tikz}

\usepackage{caption}
\usepackage{graphicx}
\usepackage{graphics}

\newtheorem{theorem}{Theorem}

\newcommand{\sign}{\operatorname{sign}}

\begin{document}

\title[]{Failure of Orthogonality\\ of Rounded Fourier Bases}

\author[]{Fran\c{c}ois Cl\'ement}
\author[]{Stefan Steinerberger}
\address{Department of Mathematics, University of Washington, Seattle}
\email{fclement@uw.edu} 
\email{steinerb@uw.edu}

\begin{abstract}  The purpose of this note is to prove estimates for
$$ \left| \sum_{k=1}^{n} \sign \left( \cos \left( \frac{2\pi a}{n}  k \right) \right) \sign \left( \cos \left( \frac{2\pi b}{n}  k \right) \right)\right|,$$
when $n$ is prime and $a,b \in \mathbb{N}$. We show that the expression can only be large if $a^{-1}b \in \mathbb{F}_n$ (or a small multiple thereof) is close to $1$. This explains some of the surprising line patterns in $A^T A$ when $A \in \mathbb{R}^{n \times n}$ is the signed discrete cosine transform. Similar results seem to exist at a great level of generality.
\end{abstract}

\maketitle

\vspace{-0pt}

\section{Introduction and Results}
\subsection{Motivation}
This paper is motivated by a completely different problem and a surprising picture. Dong and Rudelson \cite{dong} recently proved the existence of a universal constant $c>0$ such that for all $n \in \mathbb{N}$ there exists $A \in \left\{-1,1\right\}^{n \times n}$ such that $A$ is very nearly a rescaled orthogonal matrix
$$ \forall ~x \in \mathbb{R}^n \qquad \frac{1}{c}\|x\| \leq \frac{1}{\sqrt{n}} \|Ax\| \leq c \|x\|.$$
Dong-Rudelson \cite{dong} use Vinogradov's theorem that every sufficiently large odd number is the sum of three primes. A different proof \cite{stein} uses the existence of flat Littlewood polynomials \cite{little}. Neither of these approaches is elementary. 

\begin{center}
    \begin{figure}[h!]
\begin{tikzpicture}
    \node at (0,0) {\includegraphics[width=0.45\textwidth]{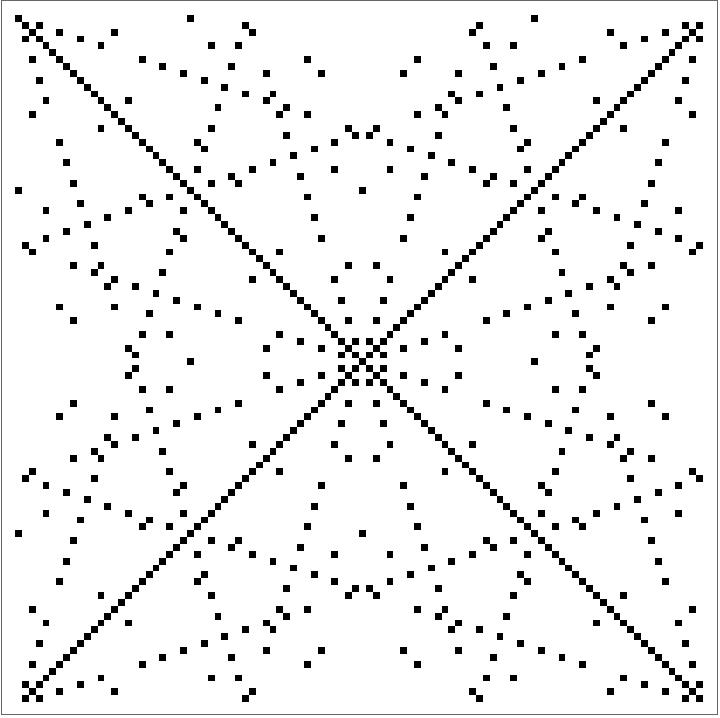}};
    \node at (6,0) {\includegraphics[width=0.45\textwidth]{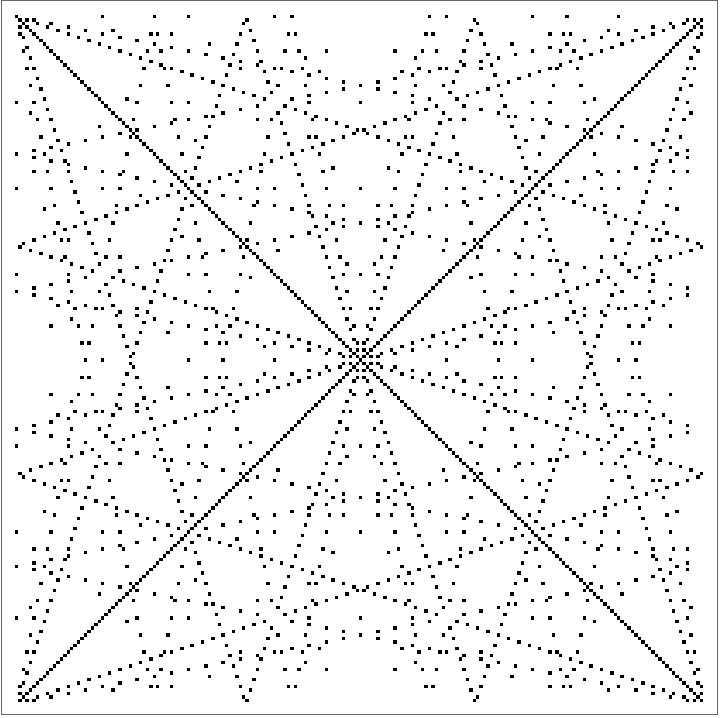}};
\end{tikzpicture}
\caption{Two motivating pictures. Left shows large entries of $A^T A$ when $n=100$, right when $n=200$. Lines start to emerge.}
    \end{figure}
\end{center}

Toying with the idea of taking a suitable orthogonal matrix $Q \in \mathbb{R}^{n \times n}$ and rounding the entries to $\left\{-1,1\right\}$  led to the discrete cosine transform and the rounded matrix
$$ A = \left( \mbox{sign}\left( \cos\left( \frac{ 2\pi}{ n} i j \right)\right)\right)_{i,j=1}^{n}.$$ 
When considering whether $\|Ax\|$ is close to an isometry, 
$\| Ax\|^2 = \left\langle Ax, Ax\right\rangle = \left\langle x, A^T A x \right\rangle$ naturally leads to $A^T A$.
Here, we encountered an amusing surprise: while $A$ does not seem to be a good approximate Hadamard matrix, the large entries of $A^T A$ seem to form a pattern that appears to `converge' as $n$ becomes large (see Fig 1). As $n$ gets larger, the `relatively large' entries of $A^T A$ appear to be concentrated on what we perceive to be `lines'. The purpose of this short note is to explain the presence of these lines and to indicate a new type of phenomenon that seems to naturally arise at a greater level of generality (see \S 1.4).

\subsection{Main Result.}
Writing out the entries of $A^T A$ leads to the expression
$$ \Sigma_{a,b} = \left| \sum_{k=1}^{n} \sign \left( \cos \left( \frac{2\pi}{n} a k \right) \right) \sign \left( \cos \left( \frac{2\pi}{n} b k \right) \right) \right|,$$
where $\mbox{sign}(x)$ is $+1$ if $x \geq 0$ and $-1$ otherwise (whether to define $\mbox{sign}(0)$ as $-1$ or $+1$ will not matter in what follows).  We will prove three types of statements that complement each other.
For the purpose of exposition, we introduce, for $x \in \mathbb{F}_n$, the abbreviation
$ \|x\| = \min\left\{x, n-x\right\}$
which measures distance from $0 \in \mathbb{F}_n$ in the toroidal sense. We first give three bounds that become increasingly complicated, a more detailed explanation is provided below.
\begin{theorem}
    Let $n$ be prime and $0 \neq a,b \in \mathbb{F}_n$. Then
    $$\Sigma_{a,b} = \left| \sum_{k=1}^{n} \sign \left( \cos \left( \frac{2\pi}{n} a k \right) \right) \sign \left( \cos \left( \frac{2\pi}{n} b k \right) \right) \right|$$
    satisfies the following estimates.
    \begin{enumerate}
        \item There is a nearly exact formula when $\|a^{-1}b\|$ is small
        $$ \Sigma_{a,b} =   \frac{n}{\|a^{-1}b\|}1 _{\left[\|a^{-1}b\|~\emph{\tiny is odd}\right]} + \mathcal{O}(\|a^{-1}b\|).$$
        \item For any $1 \leq m < n$ we have 
        $$ \Sigma_{a,b} \leq \left( \frac{2}{m+2} +\sum_{\ell=1}^{m} \frac{12}{\ell \| \ell a^{-1}b\|} \right)n + \mathcal{O}(1).$$
        \item  More generally, we have, for a suitable measure $\mu$,
     $$     \Sigma_{a,b}= \left| \frac{4n}{\pi} \sum_{\ell =1 \atop \ell ~\mbox{\tiny is odd}}^{\infty} \frac{(-1)^{(\ell-1)/2}}{\ell} \cdot \Re \widehat{\mu}(\ell) \right| + \mathcal{O}(1).$$
    
    \end{enumerate}
\end{theorem}

The first bound is easy to explain: if $1 \leq \|a^{-1} b\| \ll n,$
then we have a very precise formula: $\Sigma_{a,b}$ is close to 0 whenever $\|a^{-1} b\|$ is even and it is close to $n/\|a^{-1} b\|$ whenever $\|a^{-1} b\|$ is odd. This asymptotic also explains the lines that we observe in the picture: the two diagonals correspond to $a = b$, where $a^{-1}b = 1$ and $a = -b$ corresponding to $a^{-1}b = -1$. The next layer of lines arises from $a^{-1} b \in \left\{3, -3\right\}$ and this pattern continues. This simple bound has an error term that deteriorates when $\|a^{-1} b\|$ becomes large: this is necessary, there are examples where both $\Sigma_{a,b}$ and $\|a^{-1}b \|$ are large. The second bound completes the picture: if $\|a^{-1}b\|$ is large, the first bound gets less accurate. The second bound shows that if $\|a^{-1}b\|$ is large and if its small multiples are also large, then $\Sigma_{a,b}$ will be relatively small. One could wonder whether there is a converse result. Does the existence of a small $\ell$ with $\|\ell a^{-1}b\|$ small necessarily imply that the sum is large? This is, in some sense, answered by the third result: the truth is complicated and given by an oscillatory sum. If $\|\ell a^{-1}b\|$ is large for small $\ell$, then all the terms in the sum are small and $\Sigma_{a,b}$ is small; the converse is not necessarily true: if $\|\ell a^{-1}b\|$ is small, then $\Re \widehat{\mu}(\ell)$ may be small or may be large. Moreover, there might also be other large terms that lead to an additional cancellation. The sum in the third bound does not converge absolutely; using a standard convolution with a mollifier regularizes the sum at the cost of increasing the error term slightly, we leave the details to the reader.

\subsection{The case where $n$ is composite.}
The case where $n$ is not prime adds an additional layer of complexity. Much of the previous result remains unchanged but it is less clear whether a similarly concise description exists. One observes the emergence of additional lines on top of the ones described above (see Fig. 2). We will prove the existence of some lines of this type by an explicit method.

\begin{center}
    \begin{figure}[h!]
\begin{tikzpicture}
    \node at (0,0) {\includegraphics[width=0.45\textwidth]{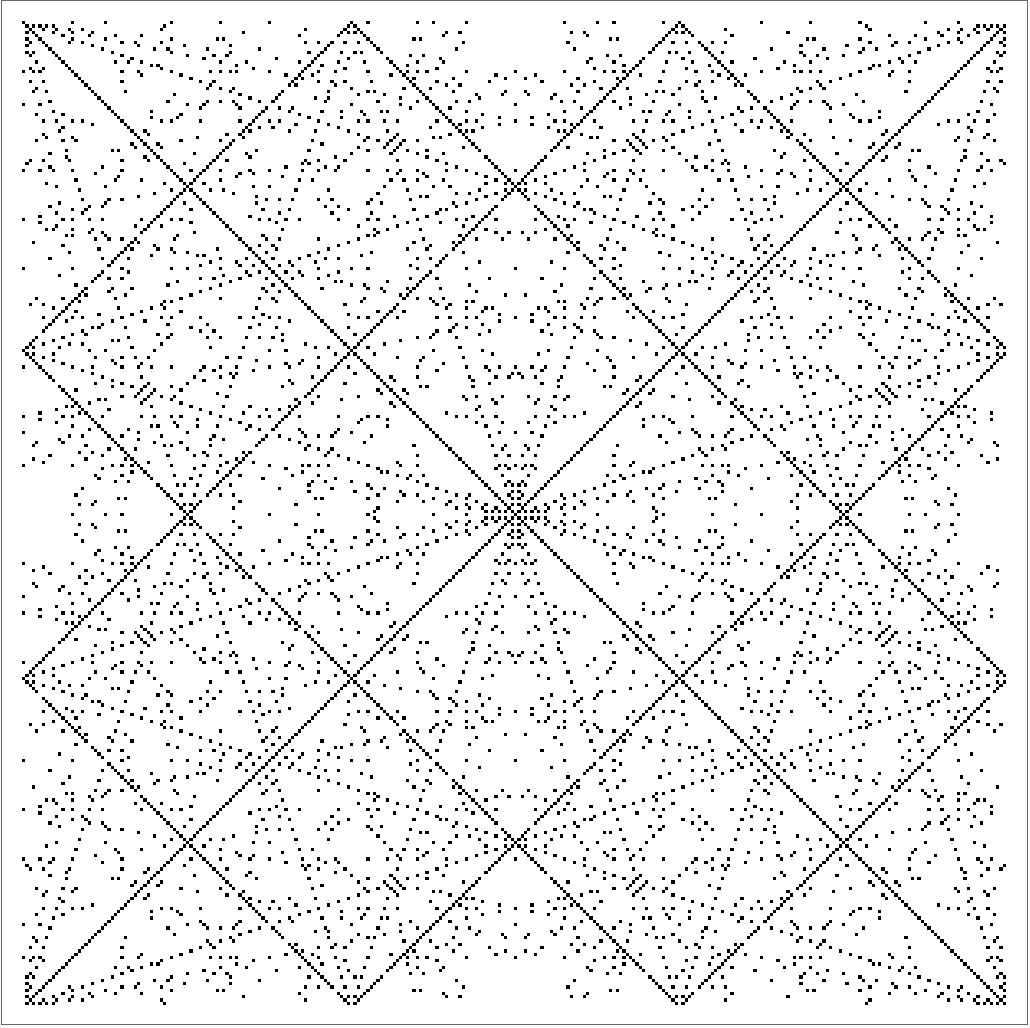}};
    \node at (6,0) {\includegraphics[width=0.45\textwidth]{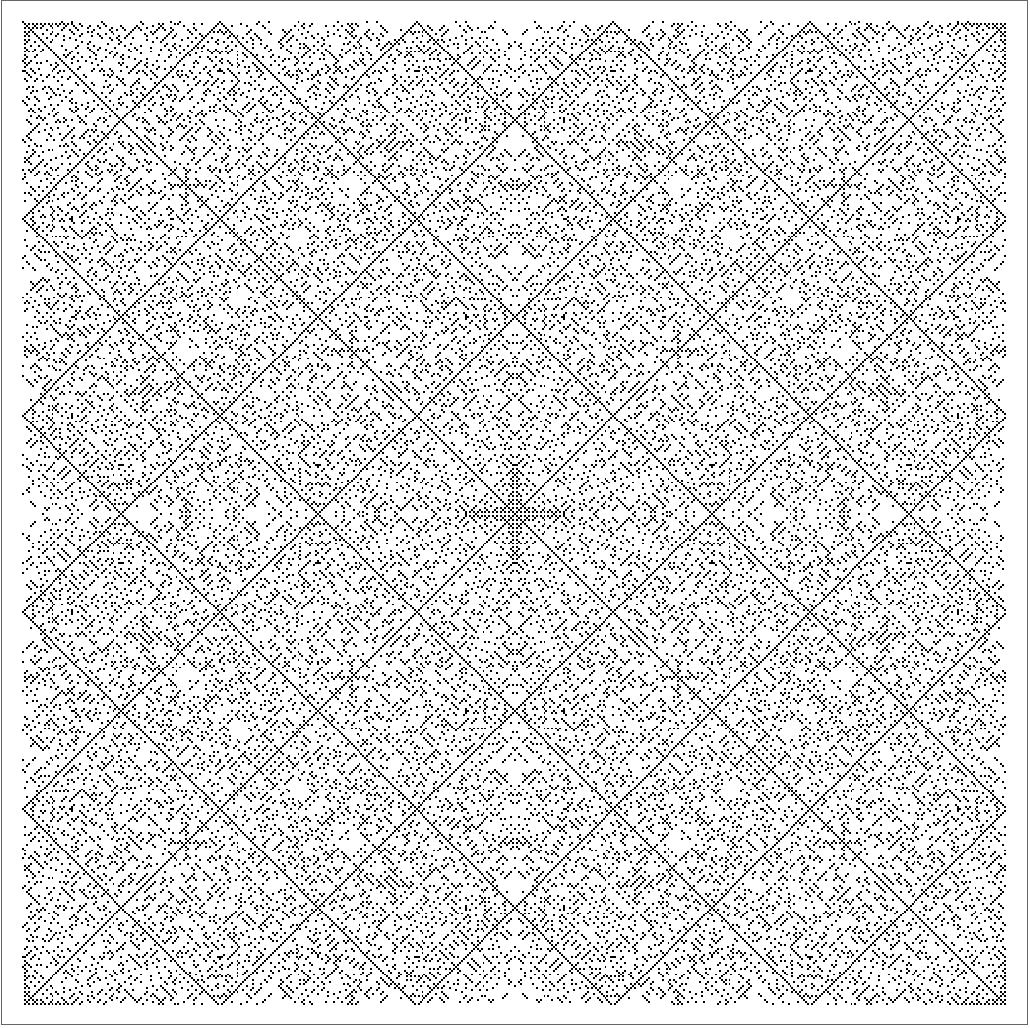}};
\end{tikzpicture}
\caption{Large entries of $A^T A$ when $n=300$ (left) and when $n=500$ (right). Additional lines emerge.}
    \end{figure}
\end{center}

\begin{theorem}
    Let $n\in \mathbb{N}$ be not prime, let $p|n$ be a prime divisor, and let $\emph{gcd}(a,n)=1$. Then for any $0<|c| \leq p-1$ 
    $$\left| \sum_{j=1}^{n} \sign\left(\cos\left( \frac{2\pi a j}{n} \right)\right)\sign\left(\cos\left( \frac{2\pi}{n}\left(a+\frac{cn}{p}\right)j\right)\right)\right| = \frac{n}{p^2}+\mathcal{O}(p).$$
\end{theorem}

We note that when $p \gtrsim n^{1/3}$, the error term starts to dominate. This result describes what happens to $\Sigma_{a,b}$ when $b-a$ is a multiple of $n/p$.

\subsection{The big picture.}
It appears that all of this is part of a much more general phenomenon.  Oscillatory orthogonal functions are sometimes not too badly approximated by their sign functions and it is conceivable that their signed version may inherit \textit{some} of the orthogonality structure (think of the Rademacher system or the Walsh system). This phenomenon appears to be quite general with our result being perhaps the simplest possible instance of it. We conclude with two more purely empirical/numerical examples.\\

\textit{1. Legendre polynomials.} The sequence of Legendre polynomials $p_n:\mathbb{R} \rightarrow \mathbb{R}$ has many different descriptions, a particularly simple one is Rodrigues' formula
$$ p_n(x) = \frac{1}{2^n n!} \frac{d^n}{dx^n} \left(x^2 - 1\right)^n.$$
These polynomials are mutually orthogonal on $[-1,1]$ with
$$ \int_{-1}^{1} p_m(x) p_n(x) dx = \frac{2}{2n+1} \delta_{mn}.$$
One could turn them orthonormal by changing the normalization constant but this does not affect anything we do here (since we take signs). One could now naturally wonder whether a quantity like
$$ Q_{mn} = \left|  \int_{-1}^{1} \mbox{sign}(p_m(x)) \mbox{sign}(p_n(x)) dx \right| $$
inherits any of the structure.

\begin{center}
    \begin{figure}[h!]
\begin{tikzpicture}
    \node at (0,0) {\includegraphics[width=0.4\textwidth]{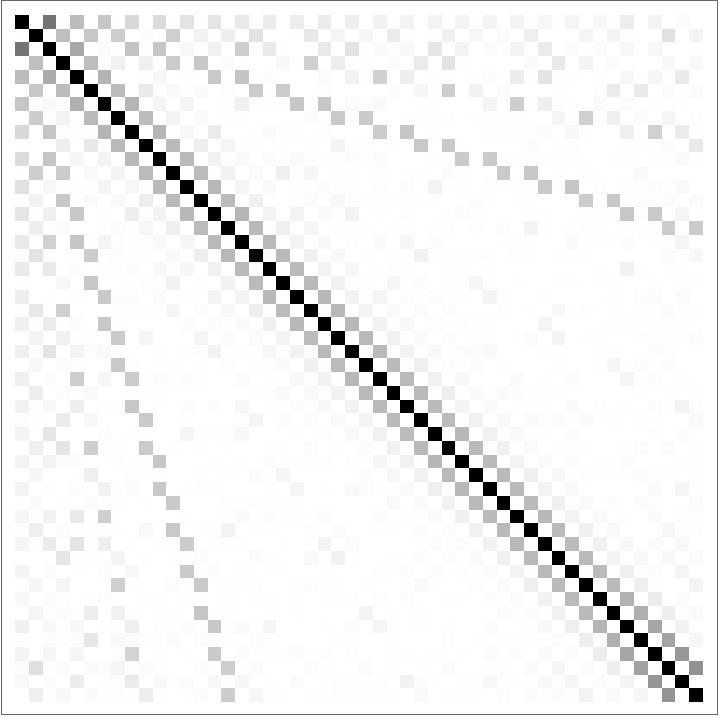}};
    \node at (6,0) {\includegraphics[width=0.4\textwidth]{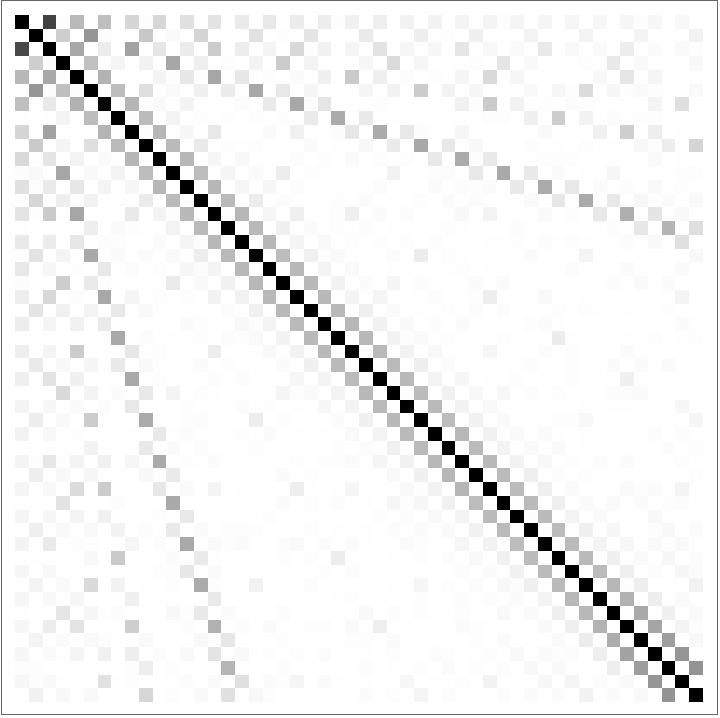}};
\end{tikzpicture}
\caption{Entries of $Q_{mn}$ (left) and $R_{mn}$ (right): larger entries are darker. Again we observe an emergence of lines.}
    \end{figure}
\end{center}

\textit{2. Chebychev polynomials.} Chebychev polynomials of the first kind are defined by
$$ T_n(\cos{n\theta}) = \cos(n \theta).$$
For $m,n \geq 0$, they are orthogonal on $[-1,1]$ with respect to the measure $(1-x^2)^{-1/2}$
$$ \int_{-1}^{1} T_n(x) T_m(x) \frac{dx}{\sqrt{1-x^2}} =  1_{m \neq n} \pi.$$
One could wonder whether 
$$ R_{mn} = \left| \int_{-1}^{1} \mbox{sign}(T_n(x)) \mbox{sign}(T_m(x)) \frac{dx}{\sqrt{1-x^2}} \right|$$
has any interesting structure. This seems to be the case. It is perhaps less surprising than one might expect: orthogonal polynomials admit a WKB asymptotic which admits a locally trigonometric approximation like the ones studied above -- though it is perhaps less clear how easy it is to make this precise. It appears that there are several interesting phenomena worthy of being further explored.

\section{Proof of Theorem 1}
\begin{proof} We ignore the last summand and only sum to $n-1$. This leads to an error of size $1$. 
We think of the set $\left\{1,2,\dots, n-1\right\} \subset \mathbb{F}_n$ as a subset of a finite field.
    The first step consists of a change of the order of summation because $k \rightarrow a^{-1}k$ is a bijection on $\left\{1,2,\dots, n-1\right\}$. The sum can be written as
\begin{align*}
\Sigma_{a,b} &= \left| \sum_{k=1}^{n-1} \sign \left( \cos \left( \frac{2\pi}{n} a k \right) \right) \sign \left( \cos \left( \frac{2\pi}{n} b k \right) \right) \right| + \mathcal{O}(1)\\
&=  \left| \sum_{k=1}^{n-1} \sign \left( \cos \left( \frac{2\pi}{n}  k \right) \right) \sign \left( \cos \left( \frac{2\pi}{n} a^{-1} b k \right) \right) \right| + \mathcal{O}(1).
\end{align*} 
We further simplify this expression by decomposing $k$ into the set where the first term is positive and its complement where the first term is negative
\begin{align*}
P &= \left\{1 \leq k \leq n-1:  \sign \left( \cos \left( \frac{2\pi}{n}  k \right) \right)  = 1\right\}\\
N &= \left\{1 \leq k \leq n-1:  \sign \left( \cos \left( \frac{2\pi}{n}  k \right) \right)  = -1\right\}.
\end{align*}
$P$ and $N$ are easy to describe: $P$ is the union of two intervals (or a single interval if we think toroidally), $N$ is an interval. 
Thus
\begin{align*}
    \Sigma_{a,b} &=  \left| \sum_{k \in P}^{}  \sign \left( \cos \left( \frac{2\pi}{n} a^{-1} b k \right) \right)  - \sum_{k \in N}^{}  \sign \left( \cos \left( \frac{2\pi}{n} a^{-1} b k \right) \right) \right| + \mathcal{O}(1).
 \end{align*}
Since
   $$ \left| \sum_{k=1}^{n-1}  \sign \left( \cos \left( \frac{2\pi}{n} a^{-1} b k \right) \right) \right| \leq 2,$$
if we know either the sum over $k \in P$ or the sum over $k \in N$, then we know both up to a $\mathcal{O}(1)$ error. Introducing
$$ A = \# \left\{  k \in N: a^{-1} b k \in N \right\},$$
we deduce
$$ \sum_{k \in N}^{}  \sign \left( \cos \left( \frac{2\pi}{n} a^{-1} b k \right) \right) =  \# N - 2A.$$
Using this equation in combination with the fact that the total sum is close to 0, 
\begin{align*} \Sigma_{a,b} &= \left| (2A - \# N) - ( \# N - 2A) \right| + \mathcal{O}(1) = \left|4A - n\right|  + \mathcal{O}(1).
\end{align*}
It remains to show that $A$ is close to $n/4$ to deduce good estimates on $\Sigma_{a,b}$. \\

\textbf{First estimate.} We start by proving our first estimate which is reasonable if and only if $\| a^{-1} b\|$ is small where, we recall, $\|x\| = \min\left\{x, n-x\right\}$. We explain the argument for $0 < a^{-1} b < n/2$, the other case is identical.

\begin{center}
    \begin{figure}[h!]
        \begin{tikzpicture}
            \draw [thick] (0,0) -- (3,0) -- (3,3) -- (0,3) -- (0,0);
            \draw [ultra thick] (0.75, -0.1) -- (0.75, 0.1);
              \draw [ultra thick] (2.25, -0.1) -- (2.25, 0.1);
              \draw [ultra thick] (-0.1, 0.75) -- (0.1, 0.75);
                \draw [ultra thick] (-0.1, 2.25) -- (0.1, 2.25);   
             \draw [dashed] (0.75,0) -- (0.75, 3);   
            \draw [dashed] (2.25,0) -- (2.25, 3);   
           \draw [dashed] (0,0.75) -- (3, 0.75);
            \draw [dashed] (0,2.25) -- (3, 2.25);   
          \node at (1.5, -0.25) {$N$};  
           \node at (0.4, -0.25) {$P$};  
           \node at (2.6, -0.25) {$P$};   
      \node at (-0.25,2.6) {$P$};   
     \node at (-0.25,0.4) {$P$};  
      \node at (-0.25,1.5) {$N$}; 
             \node at (0.4, 2.6) {+1}; 
           \node at (0.4, 1.5) {-1};  
                 \node at (0.4, 0.4) {+1}; 
       \node at (1.5, 2.6) {-1}; 
           \node at (1.5, 1.5) {+1};  
                 \node at (1.5, 0.4) {-1}; 
       \node at (2.6, 2.6) {+1}; 
           \node at (2.6, 1.5) {-1};  
                 \node at (2.6, 0.4) {+1}; 
         \draw (0,0) -- (1,3);   
         \draw (1,0) -- (2,3); 
           \draw (2,0) -- (3,3); 
        \end{tikzpicture}
        \caption{Interpreting $\mathbb{F}_n^2$ as $[0,1]^2$ and the sum as an approximation of a line integral over a closed line starting in the origin, growing with derivative 3 and looping around 3 times.}
    \end{figure}
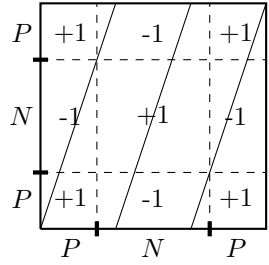
\end{center}
We are motivated by the following heuristic: if one considers $\mathbb{T}^2 \cong [0,1]^2$ and 
$$ \Phi(x,y) = \mbox{sign}(\cos(2 \pi x)) \cdot  \mbox{sign}(\cos(2 \pi y)),$$
then one would expect, for $\|a^{-1} b\|$ small, that
$$ \sum_{k=1}^{n-1} \sign \left( \cos \left( \frac{2\pi}{n}  k \right) \right) \sign \left( \cos \left( \frac{2\pi}{n} a^{-1} b k \right) \right) \sim n \int_{0}^{1} \Phi(t, \|a^{-1} b\| t) dt.$$
A computation, carried out in \cite[Lemma 3]{gon}  (see also \cite{kevin}), shows that the integral is either 0 or the inverse slope (depending on whether the slope is even or odd) 
$$ \int_{0}^{1} \Phi(t, \|a^{-1} b\| t) dt =  \frac{1}{\|a^{-1}b\|} \cdot 1_{\|a^{-1}b\|~\mbox{is odd}}.$$
It remains to quantify the error in the approximation.  The function $\Phi(x,y)$ is piecewise constant and the sampling happens at regular intervals on the $x-$axis. The only error comes from sampling closed to dashed lines: a trivial bound is thus given by understanding how often we cross one of the dashed lines and a simple estimate shows that this happens no more than $3 \|a^{-1}b\|$ times. Thus
$$ \Sigma_{a,b} =   \frac{n}{\|a^{-1}b\|}1 _{\left[\|a^{-1}b\|~\mbox{\tiny is odd}\right]} + \mathcal{O}(\|a^{-1}b\|).$$

\textit{Remark.} This argument is clearly pessimistic. It is not inconceivable that a more refined analysis might be able to improve the error estimate $\mathcal{O}(\|a^{-1}b\|)$.\\

\textbf{Second estimate.}
The second estimate uses 
\begin{align*} 
\Sigma_{a,b} = \left|4A - n\right|  + \mathcal{O}(1)
\end{align*}
together with an estimate on 
$$ A = \# \left\{  k \in N: a^{-1} b k \in N \right\} =  \# \left\{  k \in N:  \frac{1}{4} \leq \left\{\frac{a^{-1} b k}{n}\right\} \leq \frac{3}{4} \right\} + \mathcal{O}(1),$$
where $\left\{ \cdot \right\}$ is the fractional part. The problem is thus to understand the size of $$(\alpha k ~ \mbox{mod}~1)_{k=n/4}^{3n/4} \cap \left[\frac{1}{4}, \frac{3}{4} \right],$$ 
where  $ \alpha = a^{-1} b/n$. We introduce the probability measure 
$$ \mu = \frac{2}{n-1} \sum_{k=1}^{(n-1)/2} \delta_{\alpha k ~\mbox{\tiny mod}~1} $$
corresponding to the set of $(\alpha k ~ \mbox{mod}~1)_{k=n/4}^{3n/4}$ shifted by a fixed constant.
We are now interested in estimating the measure that $\mu$ assigns to an interval of length 1/2. Our argument requires an estimate for a specific interval but the subsequent argument applies uniformly to all intervals of length $1/2$.
Using Fourier series, one can write the measure as
$$ \mu = \sum_{\ell \in \mathbb{Z}} \widehat{\mu}(\ell) e^{2 \pi i \ell x} \qquad \mbox{where} \qquad \widehat{\mu}(\ell) = \frac{2}{n-1}\sum_{k=1}^{(n-1)/2} e^{-2 \pi i \ell \alpha k}.$$
The Fourier coefficient $\widehat{\mu}(\ell)$ has an unexpectedly simple form and can be computed in closed form since it is merely a geometric series
$$ \widehat{\mu}(\ell) =  \frac{2}{n-1}\sum_{k=1}^{(n-1)/2} e^{-2 \pi i \ell \alpha k} = \frac{2e^{-2 \pi i \ell \alpha}}{n-1} \frac{1 - e^{-\pi i (n-1) \ell \alpha}}{1 - e^{-2 \pi i \ell \alpha}}.$$
It satisfies the estimate
$$  \left| \widehat{\mu}(\ell) \right| = \frac{2}{n-1} \frac{\left| 1 - e^{-\pi i (n-1) \ell \alpha} \right|}{\left| 1 - e^{-2 \pi i \ell \alpha} \right|} \leq \frac{5}{n} \frac{1}{\left| 1 - e^{-2 \pi i \ell \alpha} \right|}.$$
We see that the quality of this estimate will depend on how close $\ell \alpha$ is to an integer. 
We observe that, since $ \ell \alpha = a^{-1} b \ell/n$, this is not going to be an integer for all $1 \leq \ell < n$.
Using the elementary inequality
$$ \forall~0 < x < 1 \qquad \frac{5}{|1 - e^{-2 \pi i x}|} \leq \frac{2}{\min\left\{x, 1-x\right\}},$$
we arrive at
$$  \left| \widehat{\mu}(\ell) \right| \leq \frac{2}{n}\frac{1}{\min \left\{ \ell \alpha - \left\lfloor \ell \alpha \right\rfloor , \left\lfloor \ell \alpha + 1\right\rfloor - \ell \alpha \right\}} .$$
We observe that $\ell \alpha = \ell a^{-1} b/n$ and that its distance to the nearest integer is exactly $\| \ell a^{-1} b\|/n$. Thus,
$$  \left| \widehat{\mu}(\ell) \right| \leq \frac{2}{\| \ell a^{-1} b \|}.$$
At this point, we use the Erd\H{o}s-Turan inequality in the version stated in the book of Montgomery \cite{mont}: for any $1 \leq m < n$ and any interval $J \subset \mathbb{T}$ of length $1/2$, 
$$ \left| \mu(J) - \frac{1}{2} \right| \leq   \frac{1}{m+1} + \frac{3}{n}\sum_{\ell=1}^{m} \frac{\left|\widehat{\mu}(\ell)\right|}{\ell} \leq \frac{1}{m+1} + \frac{1}{n}\sum_{\ell=1}^{m} \frac{6}{\ell \| \ell a^{-1} b\|}.$$
Using  $A =  (n-1)/2 \cdot \mu{(J)}$ for some interval $J$ of length 1/2, for all $1 \leq m < n$,
$$ \left| A - \frac{n-1}{4} \right| \leq \frac{n}{2m + 2} + \sum_{\ell=1}^{m} \frac{3}{\ell \| \ell a^{-1} b\|}.$$
Finally, invoking $\Sigma_{a,b} = |4A - n|$, one arrives at
$$ \Sigma_{a,b} \leq \frac{2n}{m+2} +\sum_{\ell=1}^{m} \frac{12}{\ell \| \ell a^{-1} b\|} + \mathcal{O}(1).$$

\textbf{Third estimate.}
The third estimate uses the fact that
  $$A =  \frac{n-1}{2} \cdot \mu\left(\left[\frac14,\frac34\right]\right) \quad \mbox{where} \quad  \mu =\sum_{\ell \in \mathbb{Z}} e^{2\pi i \ell c}\left( \frac{2}{n-1}\sum_{k=1}^{(n-1)/2} e^{-2 \pi i \ell \alpha k} \right) e^{2 \pi i \ell x},$$
  where $c \in \mathbb{T}$ is a constant coming from a global shift.
The main contribution comes from frequency $0$. We write, abbreviating the coefficients for the sake of brevity, 
$$ \mu = 1 + \sum_{\ell \in \mathbb{Z} \atop \ell \neq 0} \widehat{\mu}(\ell) e^{2 \pi i \ell x}.$$
This implies, by direct integration, that
\begin{align*}
     \mu\left(\left[\frac14,\frac34\right]\right) &= \frac{1}{2} + \sum_{\ell \in \mathbb{Z} \atop \ell \neq 0} \widehat{\mu}(\ell) \int_{1/4}^{3/4} e^{2 \pi i \ell x} dx. 
\end{align*}
We have
$$ \int_{1/4}^{3/4} e^{2 \pi i \ell x} dx  = \frac{(-1)^{ (\ell-1)/2}}{\ell \pi} 1_{\ell~\mbox{\tiny is odd}}$$
and therefore
\begin{align*}
    \mu\left(\left[\frac14,\frac34\right]\right) = \frac{1}{2} + \sum_{\ell \in \mathbb{Z} \atop \ell ~{\tiny \mbox{is odd}}} \widehat{\mu}(\ell) \frac{(-1)^{(\ell-1)/2}}{\ell \pi}.
\end{align*}
One has
$$ \widehat{\mu}(-\ell) = \overline{\widehat{\mu}(\ell)} \quad \mbox{as well as} \quad (-1)^{(-\ell-1)/2}= - (-1)^{(\ell-1)/2}$$
from which one deduces that
\begin{align*}
\widehat{\mu}(\ell) \frac{(-1)^{(\ell-1)/2}}{\ell \pi} +  \widehat{\mu}(-\ell) \frac{(-1)^{(-\ell-1)/2}}{-\ell \pi} &=     \frac{(-1)^{(\ell-1)/2}}{\ell \pi} \left( \widehat{\mu}(\ell) + \widehat{\mu}(-\ell) \right) \\
&=    \frac{2}{\pi}\frac{(-1)^{(\ell-1)/2}}{\ell} \cdot \Re \widehat{\mu}(\ell) 
\end{align*}
Summing over $\ell$ and $-\ell$ simultaneously, we get
$$  \mu\left( \left[ \frac{1}{4},\frac{3}{4}\right]\right) = \frac{1}{2} - \frac{2}{\pi}\sum_{\ell =1 \atop \ell ~\mbox{\tiny is odd}}^{\infty} \frac{(-1)^{(\ell-1)/2}}{\ell} \cdot \Re \widehat{\mu}(\ell).$$
Recalling that
  $$A =  \frac{n-1}{2} \cdot \mu\left(\left[\frac14,\frac34\right]\right) \qquad \mbox{and} \qquad  \Sigma_{a,b} = \left|4A - n\right|  + \mathcal{O}(1),$$
we have
\begin{align*}
    \Sigma_{a,b} = \left|4A - n\right|  + \mathcal{O}(1) = \left| \frac{4n}{\pi} \sum_{\ell =1 \atop \ell ~\mbox{\tiny is odd}}^{\infty} \frac{(-1)^{(\ell-1)/2}}{\ell} \cdot \Re \widehat{\mu}(\ell) \right| + \mathcal{O}(1).
\end{align*}
\end{proof}

\section{Proof of Theorem 2}
Unlike in the case where $n$ is prime, it is now possible for some cosine values to be equal to 0. In the proof, we divide the terms of the sums above based on the value of $\{cj/p\}$. For each of these values, we show that the associated cosine arguments are equally spaced in $[0,2\pi)$. The problem then becomes one of finding how many points lie in given open intervals. The cosine value equal to 0 will correspond exactly to a point being on the edge of these intervals. As such, the choice of $\sign(0)$ has little impact: we can simply count it in the counting error, as if the point were farther from the interval. These cases will always be included in the $\mathcal{O}(1)$ error.
\begin{proof}
We have
    \begin{align*}
        &\sum_{j=1}^{n} \sign\left(\cos\left( \frac{2\pi aj}{n}\right)\right)\sign\left(\cos\left(\frac{2\pi(a+cn/p)j}{n}\right)\right)\\
        &=\sum_{j=1}^n\sign\left(\cos\left(\frac{2\pi aj}{n}\right)\right)\sign\left(\cos\left(\frac{2\pi aj}{n}+2\pi\left\{\frac{cj}{p}\right\}\right)\right).
    \end{align*}
    
    We abbreviate this sum with $S$. Since $a$ and $n$ are coprime, $2\pi ja/n$ will take $n$ different values for $j \in \{1,\ldots,n\}$. 
    For $\{cj/p\} \in (0,0.5]$, the two terms $\cos\left(2\pi aj/n\right)$ and $\cos\left(2\pi\left(a+2cn/p\right)j/n\right)$  have the same sign if and only if 
    \begin{equation}\label{eq:firstint}
     \frac{2\pi aj}{n} \in \left(-\frac{\pi}{2},\frac{\pi}{2}-2\pi\left\{\frac{cj}{p}\right\}\right)\cup\left(\frac{\pi}{2},\frac{3\pi}{2}-2\pi\left\{\frac{cj}{p}\right\}\right).   
    \end{equation} 
    For $\{cj/p\} \in (0.5,1)$, 
    they have the same sign if and only if \begin{equation}\label{eq:sndint}
    \frac{2\pi aj}{n} \in \left(\frac{3\pi}{2}-2\pi\left\{\frac{cj}{p}\right\},\frac{\pi}{2}\right)\cup\left(\frac{5\pi}{2}-2\pi\left\{\frac{cj}{p}\right\},\frac{3\pi}{2}\right).
    \end{equation}
    We note that the case $cj/p=0.5$ can be included in either without loss of generality, while the case $cj/p=1$ always leads to the same signs.
The expression $\left\{cj/p\right\}$ assumes the $p$ values $0,1/p,\ldots,(p-1)/p$.
    For $i\in \{0,\ldots,p-1\}$, we introduce
    $$S_i:=\left\{j:1\leq j\leq n \quad \mbox{and} \quad cj\equiv i ~\mbox{mod}~p\right\}$$
   allowing us to change the order of summation
    $$S=\sum_{i=0}^{p-1}\sum_{j \in S_i}\sign\left(\cos\left(\frac{2\pi aj}{n}\right)\right)\sign\left(\cos\left(\frac{2\pi aj}{n}+\frac{2\pi i}{p}\right)\right).$$
The remainder of the argument now reduces to the following statement.

\begin{quote}
    \textbf{Claim.} The quantity
 $$\Sigma_i = \sum_{j \in S_i}\sign\left(\cos\left(\frac{2\pi aj}{n}\right)\right)\sign\left(\cos\left(\frac{2\pi aj}{n}+\frac{2\pi i}{p}\right)\right) 
 $$
 satisfies
$$ \Sigma_i = \begin{cases} \left(1-\frac{4i}{p}\right) \frac{n}{p} + \mathcal{O}(1) \quad &\mbox{if}~0 \leq i \leq p/2 \\
  \left(\frac{4i}{p}-3\right) \frac{n}{p} + \mathcal{O}(1) \quad &\mbox{if}~ i > p/2.\end{cases}$$
\end{quote}
The result then follows simply by summation over $0 \leq i \leq p-1$ since
    \begin{align*}
  \sum_{i=0}^{(p-1)/2} \left(\frac{n}{p}\left(1-\frac{4i}{p}\right)+\mathcal{O}(1)\right) +\sum_{i=(p+1)/2}^{p-1} \left(\frac{n}{p}\left(\frac{4i}{p}-3\right)+\mathcal{O}(1)\right) 
        =\frac{n}{p^2}+\mathcal{O}(p).
    \end{align*}
It remains to establish the claim which constitutes the remainder of the proof.\\

\textit{Proof of the Claim.} We decompose the argument into several different steps.\\
\textit{Step 1.} Since $c$ and $p$ are coprime and $p$ divides $n$, each set $S_i$ has the same cardinality ($\# S_i = n/p$). Moreover, if $j \in S_i$, then $j+p \in S_i$ and no number in between will be in $S_i$, these sets are locally simply an arithmetic progression. \\
\textit{Step 2.} The next step consists in characterizing the set of numbers that arise as arguments for the first cosine term in $\Sigma_i$, that set being $2\pi a j/n$ which leads to
$$ M_i = \left\{ \frac{aj}{n} \mod 1: j \in S_i\right\}.$$
Since the argument in the second cosine term is merely a constant shift by $2\pi i/p$ of the argument in the first term, $M_i$ describes the types of arguments that arise. We will now argue that $M_i$ is a set of $n/p$ equispaced points on $[0,1] \cong \mathbb{T}$ (when equipped with the toroidal distance). Since the $j \in S_i$ are an arithmetic progression with gap $p$, the elements of $M_i$ correspond to 
$$M_i=\left\{\frac{ad_i}{n}+\frac{ka}{n/p} \mod 1: k=0,\ldots,\frac{n}{p}-1\right\} \quad \mbox{for some} \quad d_i \in \{1,\ldots,p\}.$$
$a$ and $n/p$ are coprime and there are exactly $n/p$ elements in $M_i$. Therefore, multiplication with $a$ acts as a bijection and we can further rewrite $M_i$ as
$$M_i=\left\{\frac{ad_i}{n}+\frac{j}{n/p} \mod 1: j=0,\ldots,\frac{n}{p}-1\right\}.$$
With an initial offset of $ad_i/n$, the elements of $M_i$ can be seen to form an arithmetic progression with gap size $p/n$ on the torus.\\
\textit{Step 3.} Given an equispaced set of $m$ points on the torus $\mathbb{T} \cong [0,1]$ and an interval $J$ of length $\ell$, it is easy to see that the number of elements in $J$ is $\ell m + \mathcal{O}(1)$.\\
\textit{Step 4.} The final step now consists in finding the number of elements of $M_i$ which lead to the same sign for the second cosine. With the Equations~\eqref{eq:firstint} and~\eqref{eq:sndint}, for each $i$ we know in which intervals the elements of $M_i$ need to be for both cosines to take the same sign. By the reasoning in step 3, the number of points falling inside these intervals is given by the total number of points in $M$, $n/p$, times the length of the interval, plus a $\mathcal{O}(1)$ error term.
For $i\leq p/2$, the interval length is $(1-4i/p)$ and for $i>p/2$, the interval length is $(4i/p-3)$. For all $i$, the number of points in $M_i$ is $n/p$. As we have a union of two intervals, the error is still in $\mathcal{O}(1)$. This concludes the proof of the claim.
\end{proof}

\end{document}